\documentclass[12pt]{article}
\usepackage{amsmath,amssymb,amsbsy,amsfonts,amsthm,latexsym,
                   amsopn,amstext,amsxtra,euscript,amscd}

\newtheorem{lem}{Lemma}
\newtheorem{lemma}[lem]{Lemma}

\newtheorem{thm}{Theorem}
\newtheorem{theorem}[thm]{Theorem}

\newtheorem{ques}{Question}
\newtheorem{question}[ques]{Question}

\def\mand{\qquad\mbox{and}\qquad}

\def\\{\cr}
\def\({\left(}
\def\){\right)}
\def\[{\left[}
\def\]{\right]}
\def\<{\langle}
\def\>{\rangle}

\def\cC{{\mathcal C}}
\def\cE{{\mathcal E}}

\def\R{{\mathbb R}}
\def\eps{{\varepsilon}}

\begin{document}

\title{Integers with a large smooth divisor}

\author{
{\sc William D.~Banks} \\
{Department of Mathematics} \\
{University of Missouri} \\
{Columbia, MO 65211, USA} \\
{\tt bbanks@math.missouri.edu} \\
\and
{\sc Igor E.~Shparlinski} \\
{Department of Computing}\\
{Macquarie University} \\
{Sydney, NSW 2109, Australia} \\
{\tt igor@ics.mq.edu.au}}

\pagenumbering{arabic}

\date{}
\maketitle

\begin{abstract}
We study the function $\Theta(x,y,z)$ that counts the number of
positive integers $n\le x$ which have a divisor $d>z$ with the
property that $p\le y$ for every prime $p$ dividing $d$. We also
indicate some cryptographic applications of our results.
\end{abstract}

\section{Introduction}

For every integer $n\ge 2$, let $P^+(n)$ and $P^-(n)$ denote the
largest and the smallest prime factor of $n$, respectively, and
put $P^+(1)=1$, $P^-(1)=\infty$. For real numbers $x,y\ge 1$, let
$\Psi(x,y)$ and $\Phi(x,y)$ denote the counting functions of the
sets of \emph{$y$-smooth numbers} and \emph{$y$-rough numbers},
respectively; that is,
\begin{eqnarray*}
\Psi(x,y)&=&\#\{n\le x\ : \ P^+(n)\le y\},\\
      \Phi(x,y)&=&\#\{n\le x\ : \ P^-(n)>y\}.
\end{eqnarray*}
For a very wide range in the $xy$-plane, it is known that
$$
\Psi(x,y)\sim\varrho(u)\,x\mand\Phi(x,y)
      \sim\omega(u)\,\frac{x}{\log y},
$$
where $u$ denotes the ratio $(\log x)/\log y$, $\varrho(u)$ is the
\emph{Dickman function}, and $\omega(u)$ is the \emph{Buchstab
function}; the definitions and certain analytic properties of
$\varrho(u)$ and $\omega(u)$ are reviewed in
Sections~\ref{sec:smooth} and~\ref{sec:rough} below.

In this paper, our principal object of study is the function
$\Theta(x,y,z)$ that counts positive integers $n\le x$ for which
there exists a divisor $d\mid n$ with $d>z$ and $P^+(d)\le y$; in
other words,
$$
\Theta(x,y,z)=\#\{n\le x\ : \ n_y>z\},
$$
where $n_y$ denotes the largest $y$-smooth  
divisor  of $n$. The
function $\Theta(x,y,z)$ has been previously studied in the
literature; see~\cite{HallTen,Tenen1,Tenen2,Tenen3}.

For $x,y,z$ varying over a wide domain, we derive the first two
terms of the asymptotic expansion of $\Theta(x,y,z)$.  We show
that the main term can be naturally defined in terms of the
\emph{partial convolution} $\cC_{\omega,\varrho}(u,v)$ of
$\varrho$ with $\omega$, which is defined by
$$
\cC_{\omega,\varrho}(u,v)=\int_v^\infty\omega(u-s)\varrho(s)\,ds.
$$
Using precise estimates for $\Psi(x,y)$ and $\Phi(x,y)$, we also
identify the second term of the asymptotic expansion of
$\Theta(x,y,z)$, which is naturally expressed in terms of the
partial convolution $\cC_{\omega,\varrho'}(u,v)$ of $\varrho'$
with $\omega$:
$$
\cC_{\omega,\varrho'}(u,v)=\int_v^\infty\omega(u-s)\varrho'(s)\,ds.
$$
We use this formula to give a heuristic prediction for the density
of certain integers of cryptographic interest which appear
in~\cite{Men}.
An alternative approach, which establishes 
a two term asymptotic
formula for $\Theta(x,y,z)$ over a wider range, has been developed
recently by Tenenbaum~\cite{Tenen3}.

\begin{theorem}
\label{thm:main} For fixed $\eps>0$ and uniformly in the domain
$$
x\ge 3,\qquad y\ge\exp\{(\log\log x)^{5/3+\eps}\},
     \qquad y\log y\le z\le x/y,
$$
we have
$$
\Theta(x,y,z)=\bigl(\varrho(u)+\cC_{\omega,\varrho}(u,v)\bigr)x
-\gamma\,\cC_{\omega,\varrho'}(u,v)\,\frac{x}{\log y}
+O\bigl(\cE(x,y,z)\bigr),
$$
where $u=(\log x)/\log y$, $v=(\log z)/\log y$, $\gamma$ is the
Euler-Mascheroni constant, and
$$
\cE(x,y,z)=\frac{x}{\log y}\left\{\varrho(u-1)
     +\frac{\varrho(v)\log(v+1)}{\log y}
     +\frac{\varrho(v)}{\log(v+1)}\right\}.
$$
\end{theorem}

The proof of Theorem~\ref{thm:main} is given below in
Section~\ref{sec:proof_thm}; our principal tools are the estimates
of Lemma~\ref{lem:reciprocalsmoothX} (Section~\ref{sec:smooth})
and Lemma~\ref{lem:reciprocalroughX} (Section~\ref{sec:rough}). In
Section~\ref{sec:Crypto}, we outline some cryptographic
applications of our results.

\bigskip

\noindent{\bf Acknowledgements.} The authors would like to thank
Alfred Menezes for
bringing to our attention the cryptographic applications which
initially motivated our work.
We also thank G{\'e}rald Tenenbaum for pointing out a mistake in
the original manuscript, and for many subsequent discussions. This
work was started during a visit by W.~B.\ to Macquarie University;
the support and hospitality of this institution are gratefully
acknowledged. During the preparation of this paper, I.~S.\ was
supported in part by ARC grant DP0556431.

\section{Integers free of large prime factors}
\label{sec:smooth}

In this section, we collect various estimates for the counting
function $\Psi(x,y)$ of {\it $y$-smooth numbers\/}:
$$
\Psi(x,y)=\#\{n\le x\ : \ P^+(n)\le y\}.
$$
As usual, we denote by $\varrho(u)$ the {\it Dickman function\/};
it is continuous at $u=1$, differentiable for $u>1$, and it
satisfies the difference-differential equation
\begin{equation}
\label{eq:diffeqdickman}
      u\varrho'(u)+\varrho(u-1)=0\qquad (u>1)
\end{equation}
along with the initial condition
$$
\varrho(u)=1\qquad (0\le u\le 1).
$$
It is convenient to define $\varrho(u)=0$ for all $u<0$ so
that~\eqref{eq:diffeqdickman} is satisfied for
$u\in\R\setminus\{0,1\}$, and we also define $\varrho'(u)$ by
right-continuity at $u=0$ and $u=1$. For a discussion of the
analytic properties of $\varrho(u)$,
we refer the reader to~\cite[Chapter~III.5]{Tenen1}.

We need the following well known estimate for $\Psi(x,y)$, which
is due to Hildebrand~\cite{Hilde} (see also~\cite[Corollary~9.3,
Chapter~III.5]{Tenen1}):

\begin{lemma}
\label{lem:hildebrand} For fixed $\eps>0$ and uniformly in the
domain
$$
x\ge 3,\qquad x\ge y\ge\exp\{(\log\log x)^{5/3+\eps}\},
$$
we have
$$
\Psi(x,y)=\varrho(u)\,x\left\{1+O\(\frac{\log(u+1)}{\log
y}\)\right\},
$$
where $u=(\log x)/\log y$.
\end{lemma}

We also need the following extension of
Lemma~\ref{lem:hildebrand}, which is a special case of the results
of Saias~\cite{Saias}:

\begin{lemma}
\label{lem:saias} For fixed $\eps>0$ and uniformly in the domain
$$
x\ge 3,\qquad y\ge\exp\{(\log\log x)^{5/3+\eps}\},\qquad
     x\ge y\log y,
$$
the following estimate holds:
$$
\Psi(x,y)=\varrho(u)\,x
      +(\gamma-1)\varrho'(u)\,\frac{x}{\log y}
      +O\(\varrho''(u)\frac{x}{\log^2y}\),
$$
where $u=(\log x)/\log y$.
\end{lemma}

The following lemma provides a precise estimate for the sum
$$
S(y,z)=\sum_{\substack{d>z\\P^+(d)\le y}} \frac{1}{d}
$$
over a wide range, which is used in the proofs of
Lemmas~\ref{lem:reciprocalsmoothX} and~\ref{lem:reciprocalroughX}
below.  The sum $S(y,z)$ has been previously studied; see, for
example,~\cite{Tenen2}.

\begin{lemma} \label{lem:S(y,z)} For fixed $\eps>0$ and
uniformly in the domain
$$
y\ge 3,\qquad 1\le z\le\exp\exp\{(\log y)^{3/5-\eps}\},
$$
we have
$$
S(y,z)=\tau(v)\log y-\gamma\varrho(v)
     +O\bigl(E(y,z)\bigr),
$$
where $v=(\log z)/\log y$,
$$
\tau(v)=\int_v^\infty\varrho(s)\,ds,
$$
and
$$
   E(y,z)=\left\{
\begin{array}{ll}
        \displaystyle\frac{\varrho(v)\log(v+1)}{\log y}
         & \quad\hbox{if $z\ge y\log y$;} \\ \\
        \displaystyle z^{-1}+\frac{\log\log y}{\log y}
         & \quad\hbox{if $z<y\log y$.}
\end{array}
\right.
$$
\end{lemma}

\begin{proof}
Let $Y=y\log y$. First, suppose that $z>Y$, and put
$$
T=\frac{\exp\{(\log y)^{3/5-\eps/2}\}}{\log y}.
$$
By partial summation, it follows that
\begin{equation}
\label{eq:expand_S(y,z)}
\begin{split}
S(y,z)&=\sum_{\substack{z<d\le y^T\\P^+(d)\le y}}
\frac{1}{d}+S(y,y^T)\\
&=\frac{\Psi(y^T,y)}{y^T}-\frac{\Psi(z,y)}{z}
     +\log y\int_v^{T}\frac{\Psi(y^s,y)}{y^s}\,ds+S(y,y^T).
\end{split}
\end{equation}
By Lemma~\ref{lem:hildebrand}, we have the estimate
$$
\frac{\Psi(z,y)}{z}
   =\varrho(v)+O\(\frac{\varrho(v)\log(v+1)}{\log y}\).
$$
Also, by our choice of $T$ we have
\begin{equation}
\label{eq:S(y,z)-A} \frac{\Psi(y^T,y)}{y^T}\ll\varrho(T)
     \ll\frac{\varrho(v)\log(v+1)}{\log y}.
\end{equation}
The following bound is given in the proof
of~\cite[Corollaire~2]{Tenen2}:
$$
S(y,y^T)=\sum_{\substack{d>y^T\\P^+(d)\le y}}\frac{1}{d}\ll
\varrho(T)e^{\eps T}+y^{-(1-\eps)T},
$$
from which we deduce that
\begin{equation}
\label{eq:S(y,z)-C} S(y,y^T)\ll\frac{\varrho(v)\log(v+1)}{\log y}.
\end{equation}
To estimate the integral in~\eqref{eq:expand_S(y,z)}, we apply
Lemma~\ref{lem:saias} and write
$$
\int_v^{T}\frac{\Psi(y^s,y)}{y^s}\,ds=I_1+I_2+O(I_3),
$$
where
\begin{eqnarray*}
I_1&=&\int_v^T\varrho(s)\,ds=\tau(v)-\tau(T),\\
I_2&=&\frac{(\gamma-1)}{\log y}\int_v^T\varrho'(s)\,ds
     =\frac{(\gamma-1)(\varrho(T)-\varrho(v))}{\log y},\\
I_3&=&\frac{1}{\log^2y}\int_v^T\varrho''(s)\,ds
     =\frac{\varrho'(T)-\varrho'(v)}{\log^2y}.
\end{eqnarray*}
Since $|\varrho'(v)|\asymp\varrho(v)\log(v+1)$, and
$$
\tau(T) \ll \varrho(T) \ll \frac{\varrho(v)\log(v+1)}{\log^2 y},
$$
it follows that
\begin{equation}
\label{eq:S(y,z)-D} \int_v^{T}\frac{\Psi(y^s,y)}{y^s}\,ds
     =\tau(v)-\frac{(\gamma-1)\varrho(v)}{\log y}
     +O\(\frac{\varrho(v)\log(v+1)}{\log^2y}\).
\end{equation}
Inserting the estimates~\eqref{eq:S(y,z)-A}, \eqref{eq:S(y,z)-C}
and~\eqref{eq:S(y,z)-D} into~\eqref{eq:expand_S(y,z)}, we obtain
the desired estimate in the case $z>Y$.

Next, suppose that $y\le z\le Y$, and put
$$
V=\frac{\log Y}{\log y}=1+\frac{\log\log y}{\log y}.
$$
Since $\varrho(s)=1-\log s$ for $1\le s\le 2$, we have
$$
1 \ge \varrho(v) \ge \varrho(V)=1+O\(\frac{\log\log y}{\log y}\);
$$
therefore,
\begin{equation}
\label{eq:S(y,z)-E} \varrho(v)-\varrho(V)\ll\frac{\log\log y}{\log
y}.
\end{equation}
By partial summation, it follows that
\begin{equation}
\label{eq:S(y,z)-F}
\begin{split}
S(y,z)
     &=\sum_{\substack{z<d\le Y\\P^+(d)\le y}}
     \frac{1}{d}+S(y,Y)\\
     &=\frac{\Psi(Y,y)}{Y}-\frac{\Psi(z,y)}{z}
     +\log y\int_v^{V}\frac{\Psi(y^s,y)}{y^s}\,ds
     +S(y,Y).
\end{split}
\end{equation}
Using Lemma~\ref{lem:hildebrand} together
with~\eqref{eq:S(y,z)-E}, it follows that
\begin{equation}
\label{eq:S(y,z)-G} \frac{\Psi(Y,y)}{Y}-\frac{\Psi(z,y)}{z}
     =\varrho(V)-\varrho(v)+O\(\frac{1}{\log y}\)
     \ll\frac{\log\log y}{\log y}.
\end{equation}
Applying the estimate from the previous case, we also have
\begin{equation}
\label{eq:S(y,z)-H} S(y,Y)=\tau(V)\log Y-\gamma\varrho(V)
     +O\(\frac{1}{\log y}\).
\end{equation}
To estimate the integral in~\eqref{eq:S(y,z)-F}, we use
Lemma~\ref{lem:hildebrand} again and write
$$
\int_v^{V}\frac{\Psi(y^s,y)}{y^s}\,ds=I_4+O(I_5),
$$
where
\begin{eqnarray*}
I_4&=&
\int_v^{V}\varrho(s)\,ds=\tau(v)-\tau(V),\\
I_5&=&\frac{1}{\log y}\int_v^{V}ds
     =\frac{\log(Y/z)}{\log^2y}\ll\frac{\log\log y}{\log^2y}.
\end{eqnarray*}
Therefore,
\begin{equation}
\label{eq:S(y,z)-I} \int_v^{V}\frac{\Psi(y^s,y)}{y^s}\,ds
=\tau(v)-\tau(V)+O\(\frac{\log\log y}{\log^2y}\).
\end{equation}
Inserting the estimates~\eqref{eq:S(y,z)-G}, \eqref{eq:S(y,z)-H}
and~\eqref{eq:S(y,z)-I} into~\eqref{eq:S(y,z)-F}, and taking into
account~\eqref{eq:S(y,z)-E}, we obtain the stated estimate for
$y\le z\le Y$.

Finally, suppose that $1\le z < y$. In this case,
\begin{equation}
\label{eq:S(y,z)-J} S(y,z)=\sum_{z<d\le y}\frac{1}{d}+S(y,y).
\end{equation}
By partial summation, we have
\begin{eqnarray*}
\sum_{z<d\le y}\frac{1}{d}&=&\log y-\log z+O(z^{-1})
     =(1-v)\log y+O(z^{-1})\\
&=&\log y\int_v^1\varrho(s)\,ds+O(z^{-1}) =(\tau(v)-\tau(1))\log
y+O(z^{-1}).
\end{eqnarray*}
Applying the estimate from the previous case, we also have
$$
S(y,y)=\tau(1)\log y-\gamma\varrho(1)
     +O\(\frac{\log\log y}{\log y}\).
$$
Inserting these estimates into~\eqref{eq:S(y,z)-J}, and using the
fact that $\varrho(v)=\varrho(1)=1$, we obtain the desired result.
This completes the proof.
\end{proof}

\begin{lemma}
\label{lem:reciprocalsmoothX} For fixed $\eps>0$ and uniformly in
the domain
$$
x\ge 3,\qquad y\ge\exp\{(\log\log x)^{5/3+\eps}\},
     \qquad 1\le z\le x/y,
$$
we have
$$
\sum_{\substack{z<d\le x/y\\P^+(d)\le y}}\frac{\varrho(u-u_d)}{d}
     \ll\cC_{\varrho,\varrho}(u,v)\log(u+1)+\varrho(u-v)\varrho(v)+\varrho(u-1),
$$
where $u=(\log x)/\log y$, $v=(\log z)/\log y$, $u_d=(\log d)/\log
y$ for every integer $d$ in the sum, and
$$
\cC_{\varrho,\varrho}(u,v)=\int_v^\infty\varrho(u-s)\varrho(s)\,ds.
$$
\end{lemma}

\begin{proof}
By partial summation, we have
$$
\sum_{\substack{z<d\le x/y\\P^+(d)\le y}}\frac{\varrho(u-u_d)}{d}
=S(y,x/y)-\varrho(u-v)S(y,z)+\int_v^{u-1}\varrho'(u-s)S(y,y^s)\,ds.
$$
Lemma~\ref{lem:S(y,z)} implies that
\begin{equation*}
\begin{split}
S(y,x/y)&=\tau(u-1)\log y+O\bigl(\varrho(u-1)\bigr),\\
S(y,z)&=\tau(v)\log y+O\bigl(\varrho(v)\bigr),
\end{split}
\end{equation*}
and
$$
\int_v^{u-1}\varrho'(u-s)S(y,y^s)\,ds=I_1\log y+O(I_2),
$$
where
\begin{equation*}
\begin{split}
I_1&=\int_v^{u-1}\varrho'(u-s)\tau(s)\,ds
=\varrho(u-v)\tau(v)-\tau(u-1)+\cC_{\varrho,\varrho}(u,v),\\
I_2&=\int_v^{u-1}\bigl|\varrho'(u-s)\bigr|
     \varrho(s)\,ds.
\end{split}
\end{equation*}
Finally, using the bound
$$
\bigl|\varrho'(t)\bigr|\ll\varrho(t)\log(t+1)\qquad(t>1),
$$
we see that
$$
I_2\ll\log(u+1)\int_v^{u-1}\varrho(u-s)
\varrho(s)\,ds\le\cC_{\varrho,\varrho}(u,v)\log(u+1).
$$
Putting everything together, the result follows.
\end{proof}

\section{Integers free of small prime factors}
\label{sec:rough}

In this section, we collect various estimates for the counting
function $\Phi(x,y)$ of {\it $y$-rough numbers\/}:
$$
\Phi(x,y)=\#\{n\le x\ : \ P^-(n)>y\}.
$$
As usual, we denote by $\omega(u)$ the {\it Buchstab function\/};
for $u>1$, it is the unique continuous solution to the
difference-differential equation
\begin{equation}
\label{eq:diffeqbuch}
      \bigl(u\omega(u)\bigr)'=\omega(u-1)\qquad (u>2)
\end{equation}
with initial condition
$$
u\omega(u)=1\qquad (1\le u\le 2).
$$
It is convenient to define $\omega(u)=0$ for all $u<1$ so
that~\eqref{eq:diffeqbuch} is satisfied for
$u\in\R\setminus\{1,2\}$, and we also define $\omega'(u)$ by
right-continuity at $u=1$ and $u=2$. For a discussion of the
analytic properties of $\omega(u)$,
we refer the reader to~\cite[Chapter~III.6]{Tenen1}

The next result follows from~\cite[Corollary~7.5,
Chapter~III.6]{Tenen1}:

\begin{lemma}
\label{lem:smallprimefree} For fixed $\eps>0$ and uniformly in the
domain
$$
x\ge 3,\qquad x\ge y\ge\exp\{(\log\log x)^{5/3+\eps}\},
$$
the following estimate holds:
$$
\Phi(x,y)=\bigl(x\omega(u)-y\bigr)\frac{e^\gamma}{\zeta(1,y)}
      +O\(\frac{x\varrho(u)}{\log^2y}\),
$$
where $u=(\log x)/\log y$, and $\zeta(1,y)=\prod_{p\le
y}\(1-p^{-1}\)^{-1}$.
\end{lemma}

\begin{lemma}
\label{lem:reciprocalroughX} For fixed $\eps>0$ and uniformly in
the domain
$$
x\ge 3,\qquad y\ge\exp\{(\log\log x)^{5/3+\eps}\},
     \qquad 1\le z\le x/y,
$$
we have
$$
\sum_{\substack{z<d\le x/y\\P^+(d)\le y}}\frac{\omega(u-u_d)}{d}
     =\cC_{\omega,\varrho}(u,v)\log y
     -\gamma\,\cC_{\omega,\varrho'}(u,v)+O\bigl(E(y,z)\bigr),
$$
where $u=(\log x)/\log y$, $v=(\log z)/\log y$, $u_d=(\log d)/\log
y$ for every integer $d$ in the sum, and $E(y,z)$ is the error
term of Lemma~\ref{lem:S(y,z)}.
\end{lemma}

\begin{proof}
By partial summation, it follows that
$$
\sum_{\substack{z<d\le x/y\\P^+(d)\le y}}\frac{\omega(u-u_d)}{d}
=S(y,x/y)-\omega(u-v)S(y,z)+\int_v^{u-1}\omega'(u-s)S(y,y^s)\,ds.
$$
By Lemma~\ref{lem:S(y,z)} we have the estimates
$$
S(y,x/y)=\tau(u-1)\log y-\gamma\varrho(u-1)+O\bigl(E(y,x/y)\bigr)
$$
and
$$
S(y,z)=\tau(v)\log y-\gamma\varrho(v)+O\bigl(E(y,z)\bigr).
$$
Also,
$$
\int_v^{u-1}\omega'(u-s)S(y,y^s)\,ds=I_1\log y-\gamma I_2+O(I_3),
$$
where
\begin{equation*}
\begin{split}
I_1&=\int_v^{u-1}\omega'(u-s)\tau(s)\,ds
=\omega(u-v)\tau(v)-\tau(u-1)+\cC_{\omega,\varrho}(u,v),\\
I_2&=\int_v^{u-1}\omega'(u-s)\varrho(s)\,ds
=\omega(u-v)\varrho(v)-\varrho(u-1)+\cC_{\omega,\varrho'}(u,v),\\
I_3&=\frac{1}{\log y}\int_v^{u-1}\bigl|\omega'(u-s)\bigr|
     E(y,y^s)\,ds.
\end{split}
\end{equation*}
Putting everything together, we see that the stated estimate
follows from the bound
\begin{equation}
\label{eq:EEIErough} E(y,x/y)+\omega(u-v)E(y,z)+I_3\ll E(y,z).
\end{equation}
To prove this, observe that $E(y,z_1)\ll E(y,z_2)$ holds for all
$z_1\ge z_2\ge 1$.  Therefore, $E(y,x/y)\ll E(y,z)$, and
$$
I_3\ll\frac{E(y,z)}{\log
y}\int_v^{u-1}\bigl|\omega'(u-s)\bigr|\,ds \ll\frac{E(y,z)}{\log
y}.
$$
Taking into account the fact that $\omega(u-v)\asymp 1$, we derive
the bound~\eqref{eq:EEIErough}, and this completes the proof.
\end{proof}

\section{Proof of Theorem~\ref{thm:main}}
\label{sec:proof_thm}

For fixed $y$, every positive integer $n$ can be uniquely
decomposed as a product $n=de$, where $P^+(d)\le y$ and
$P^-(e)>y$. Therefore,
\begin{eqnarray*}
\Theta(x,y,z)&=&
      \sum_{\substack{z<d\le x\\P^+(d)\le y}}
      \sum_{\substack{e\le x/d\\P^-(e)>y}}1
      =\sum_{\substack{z<d\le x\\P^+(d)\le y}}\Phi(x/d,y)\\
      &=&\Psi(x,y)-\Psi(x/y,y)+
      \sum_{\substack{z<d\le x/y\\P^+(d)\le y}}\Phi(x/d,y).
\end{eqnarray*}
Using Lemma~\ref{lem:hildebrand}, it follows that
$$
\Psi(x,y)-\Psi(x/y,y)=\varrho(u)\,x+O\(\frac{\varrho(u-1)\,x}{\log
y}\).
$$
By Lemma~\ref{lem:smallprimefree}, we also have
\begin{equation}
\label{eq:lastexpansion}
\begin{split}
\sum_{\substack{z<d\le x/y\\P^+(d)\le y}}&\Phi(x/d,y)\\
      &=\sum_{\substack{z<d\le x/y\\P^+(d)\le y}}
      \left\{\(\frac{x\,\omega(u-u_d)}{d}-y\)\frac{e^\gamma}{\zeta(1,y)}
      +O\(\frac{x\varrho(u-u_d)}{d\log^2y}\)\right\}\\
      &=\frac{e^\gamma x}{\zeta(1,y)}
       \sum_{\substack{z<d\le x/y\\P^+(d)\le y}}\frac{\omega(u-u_d)}{d}
      -\frac{e^\gamma y}{\zeta(1,y)}\bigl\{\Psi(x/y,y)-\Psi(z,y)\bigr\}\\
      &\qquad\qquad+\quad
O\Biggl(\frac{x}{\log^2y}\sum_{\substack{z<d\le x/y\\
      P^+(d)\le y}}\frac{\varrho(u-u_d)}{d}\Biggl).
\end{split}
\end{equation}
Applying Lemma~\ref{lem:hildebrand} again, we have
$$
-\frac{e^\gamma y}{\zeta(1,y)}\bigl\{\Psi(x/y,y)-\Psi(z,y)\bigr\}
     \ll\frac{\varrho(u-1)\,x}{\log y}.
$$
Inserting the estimates of Lemmas~\ref{lem:reciprocalsmoothX}
and~\ref{lem:reciprocalroughX} into~\eqref{eq:lastexpansion}, and
making use of the trivial estimate
$$
\cC_{\varrho,\varrho}(u,v)\log(u+1)\ll\log y
\int_v^\infty\varrho(s)\,ds\ll\frac{\varrho(v)\log y}{\log(v+1)}.
$$
it is easy to see that
$$
\Theta(x,y,z)=\biggl(\varrho(u)+\cC_{\omega,\varrho}(u,v)
\frac{e^\gamma\log y}{\zeta(1,y)}\biggl)x
-\gamma\,\cC_{\omega,\varrho'}(u,v)\frac{e^\gamma x}{\zeta(1,y)}
+O\bigl(\cE(x,y,z)\bigr).
$$
To complete to proof, we use the estimate (see~\cite{Vino}):
$$
\zeta(1,y)=e^\gamma\log y\bigl(1+\exp\{-c(\log y)^{3/5}\}\bigr),
$$
which holds for some absolute constant $c>0$, together with the
trivial estimate
$$
\max\bigl\{\cC_{\omega,\varrho}(u,v),\cC_{\omega,\varrho'}(u,v)\bigr\}\ll
\int_v^\infty\varrho(s)\,ds\ll\frac{\varrho(v)}{\log(v+1)}.
$$

\section{Cryptographic applications}
\label{sec:Crypto}

Suppose that two primes $p$ and $q$ are selected for use in the
\emph{Digital Signature Algorithm} (see, for example,~\cite{MOV})
using the following standard method:
\begin{itemize}
\item Select a random $m$-bit prime $q$;

\item Randomly generate $k$-bit integers $n$ until a prime
$p=2nq+1$ is reached.
\end{itemize}
The \emph{large subgroup attack} described
in~\cite[Section~3.2.2]{Men} leads one naturally to consider the
following question: What is the probability $\eta(k,\ell, m)$
that~$n$ has a divisor $s > q$ which is $2^\ell$-smooth?

It is natural to expect that the proportion of those integers in
the set $\{2^{k-1}\le n<2^k\}$ having a large smooth divisor
should be roughly the same as the proportion of integers in
$$
\left\{2^{k-1}\le n<2^k~:~n=(p-1)/(2q)\text{~for some
prime~}p\equiv 1\pmod{2q}\right\}
$$
having a large smooth divisor.  
Accordingly, we expect that the
probability $\eta(k,\ell,m)$ is reasonably close to
$$
\frac{\Theta(2^k,2^\ell,2^m)-\Theta(2^{k-1},2^\ell,2^m)}{2^{k-1}}.
$$
Theorem~\ref{thm:main} then suggests that
$$
\eta(k,\ell,m)\approx
     2\,\wp(k,\ell,m)-\wp(k-1,\ell,m),
$$
where
$$
\wp(k,\ell,m)=\varrho(k/\ell)+\cC_{\omega,\varrho}(k/\ell,m/\ell)
-\frac{\gamma\,\cC_{\omega,\varrho'}(k/\ell,m/\ell)}{\ell\log 2}.
$$
In particular, the most interesting choice of parameters at the
present time is $k = 863$, $\ell = 80$,  and $m = 160$ (which
produces a $1024$-bit prime $p$), for which expect that
$\eta(863,80,160)\approx 0.09576$.

\end{document}